\documentstyle[amscd,amssymb,verbatim,12pt]{amsart}

\theoremstyle{plain}
\newtheorem{Thm}{Theorem}[section]
\newtheorem{Lem}[Thm]{Lemma}
\newtheorem{Def}[Thm]{Definition}
\newtheorem{Cor}[Thm]{Corollary}
\newtheorem{Prop}[Thm]{Proposition}

\theoremstyle{remark}

\newtheorem{remark}{Remark}

\newcommand{\cha}{\operatorname{char}}

\newcommand{\gal}{\operatorname{Gal}}
\newcommand{\im}{\operatorname{im}}

\begin{document}
\title[Galois groups over nonrigid fields]{Galois groups over nonrigid
fields}
\author{Wenfeng Gao}
\address{537 Bellevue Way SE, Apt. 110, Bellevue WA 98004}
\email{wgao@@sprintmail.com}
\author{David B. Leep}
\address{Department of Mathematics\\University of Kentucky\\Lexington, KY
40506-0027 U.S.A}
\email{leep@@ms.uky.edu}
\author{J\'{a}n Min\'{a}\v{c}}
\address{Department of Mathematics\\University of Western Ontario\\London,
Ontario, Canada N6A 5B7}
\email{minac@@uwo.ca}
\thanks{The third author was partially supported by NSERC and the special
Dean of Science Fund at UWO. He also gratefully acknowledges the support of
MSRI in Berkeley during the fall of 1999.}
\author{Tara L. Smith}
\address{Department of Mathematical Sciences\\University of
Cincinnati\\Cincinnati, Ohio 45221-0025 U.S.A}
\email{tsmith@@math.uc.edu}

\dedicatory{Dedicated to Paulo Ribenboim}



\begin{abstract}
Let $F$ be a field with $\cha F \neq 2$.  We show that $F$ is a nonrigid
field if and only if certain small $2$-groups occur as Galois groups over
$F$.  These results provide new \lq\lq automatic realizability" results for
Galois groups over $F$.  The groups we consider demonstrate the inequality
of two particular metabelian $2$-extensions of $F$ which are unequal
precisely when $F$ is a nonrigid field.  Using known results on connections
between rigidity and existence of certain valuations, we obtain
Galois-theoretic criteria for the existence of these valuations.
\end{abstract}

\maketitle

\section{Introduction} \label{S:intro}

Let $F$ be a field with $\cha F \neq 2$.  The goal of this paper is to
identify basic Galois groups that must occur if $F$ is not a rigid field.
If $|\dot{F}/\dot{F}^2| \geq 8$ then it is known that $F$ is not rigid if
and only if a certain group (often denoted $DC$) of order $16$ occurs as a
Galois group over $F$ ([MS]).  It has also been shown ([LS], simplifying an
earlier proof in [AGKM]) that $F$ is rigid if and only if
$F^{(3)}=F^{\{3\}}$ (all notation defined below).  In this paper we
identify two groups, denoted $G_1$ and $G_2$, where $|G_1| =32$,
$|G_2|=64$, with the property that a field $F$ is nonrigid if and only if
at least one of $G_1$ and $G_2$ is realizable as a Galois group over $F$.
In each case  the corresponding Galois extension lies in $F^{\{3\}}$ but
not in $F^{(3)}$.  Moreover, it is shown that extensions realizing $G_1$
and $G_2$ correspond, respectively, to extensions realizing the groups $D
\curlywedge C$ and $D \curlywedge D$, which lie inside $F^{(3)}$.  This
correspondence gives new nontrivial ``automatic realizability" results for
Galois groups.  We also point out how to detect nonrigid elements in $F$
using the groups $G_1$ and $G_2$.  Finally, since rigidity conditions on
field elements are known to correspond to the existence of certain
valuations ([W]), we can provide a Galois-theoretic interpretation of
valuation-theoretic results.

We begin in the following section with necessary notation and terminology,
as well as an explanation of the extensions $F^{(3)}$ and $F^{\{3\}}$ and
properties of their Galois groups.   In section 3 we examine the group $D
\curlywedge C$ and its associated group $G_1$, and in section 4 the groups $D
\curlywedge D$ and $G_2$ are considered.  Section 5 gives the main theorems
connecting nonrigidity to the realizability of the groups $D\curlywedge C$,
$G_1$, $D\curlywedge D$ and $G_2$, as well as the Galois-theoretic
interpretation of the existence of a nonrigid element.  We conclude in section
6 with a detailed analysis of the situation when $|\dot{F}/\dot{F}^2| =4$.
This is the case where the group $DC$ is not realizable as a Galois group
over $F$, and yet $F$ may be nonrigid.

\section{Notation, terminology and Preliminaries}

We let $\dot{F} = F \backslash \{0\}$, the group of nonzero elements of
$F$.  If $a_1, \ldots a_n \in \dot{F}$, the notation $\langle a_1, \ldots ,
a_n \rangle$ denotes the quadratic form $a_1x_1^2 + \cdots + a_nx_n^2$.  If
$q$ is a quadratic form defined over $F$, the notation $D_F(q)$ denotes the
set of elements in $\dot{F}$ represented by $q$ over $F$.  Basic results on
quadratic forms can be found in [La].

If $G$ is a group and $\sigma, \tau \in G$, then $[\sigma, \tau ]$ denotes
the commutator $\sigma^{-1}\tau^{-1}\sigma\tau$.

The notion of a rigid field and a rigid element are important in what
follows.  An element $a \in \dot{F} \backslash \pm F^2$ is {\em rigid} if
$D_F(\langle 1,a \rangle) = \dot{F}^2 \cup a\dot{F}^2$ and is {\em
nonrigid} otherwise.  An element $a \in \dot{F} \backslash \pm F^2$ is {\em
double rigid} if $a$ and $-a$ are both rigid.  A field $F$ is {\em rigid}
if each $a \in \dot{F} \backslash \pm F^2$ is rigid and is {\em nonrigid}
otherwise.

\begin{Lem} The following are equivalent.
\begin{enumerate}
\item $F$ is nonrigid.
\item There exists $a \in \dot{F} \backslash F^2$ such that $|D_F(\langle
1,-a \rangle)  /\dot{F}^2| \geq 4$.
\item There exists $a \in \dot{F} \backslash \pm F^2$ such that
$|D_F(\langle 1,-a \rangle)  /\dot{F}^2| \geq 4$.
\end{enumerate}
\end{Lem}

\begin{pf}
It is clear that (1) and (3) are equivalent and that (3) implies (2).
Now assume that (2) holds but that each element in $\dot{F} \backslash \pm
F^2$ is rigid.  Then we can assume $a=-1$ in (2) and that $-1 \notin F^2$.
Then there exists $b \in \dot{F} \backslash \pm F^2$ such that $b \in
D_F(\langle 1,1 \rangle)$.  Then $-1 \in D_F(\langle 1, -b  \rangle) =
\dot{F}^2 \cup -b\dot{F}^2$, which is a contradiction.
\end{pf}

An element $a \in \dot{F}$ is {\em basic} if either $a \in \pm F^2$ or $a$
is not double rigid.  The set of all basic elements in $F$ is denoted
$B(F)$.  Thus $\dot{F} = B(F) \cup \{\text{double rigid elements}\}$ is a
disjoint union and $B(F)= \dot{F}^2 \cup -\dot{F}^2$ if and only if $F$ is a
rigid field.

Let $v : F \to \Gamma \cup \{\infty\}$ be a valuation on $F$ with valuation
ring $A$, maximal ideal $M$, residue field $k=A/M$ and value group
$\Gamma$, where $\Gamma$ is an ordered abelian group.  Let $U=A \backslash
M$ be the group of units of $A$.  Then $\Gamma \cong \dot{F}/U$.  The
valuation $v$ is called $2$-henselian if $v$ has a unique extension to the
quadratic closure of $F$.  If $\cha F \neq 2$ and $1+M \subseteq
\dot{F}^2$, then $v$ is $2$-henselian ([W], Lemma 4.3).  The valuation
$v$ is called $2$-divisible if $\Gamma$ is a $2$-divisible group.  Since
$\dot{F}^2 U$ is the set of elements in $F$ whose valuation lies in
$2\Gamma$, it follows that $v$ is $2$-divisible if and only if $\dot{F}
=\dot{F}^2 U$.

The following theorem first appeared in [W], Theorem 4.4 (3).

\begin{Thm}
Suppose $F$ is nonrigid.  Then there exists a valuation $v$ on $F$ such
that $1+M \subseteq \dot{F}^2$ and $B(F) = \dot{F}^2 U$.
\end{Thm}

\begin{Cor}
Suppose $F$ is nonrigid.  The following are equivalent.
\begin{enumerate}
\item There exists a valuation $v$ on $F$ such that $1+M \subseteq
\dot{F}^2$,
$B(F) = \dot{F}^2 U$ and $v$ is not $2$-divisible.
\item $F$ contains a double rigid element.
\end{enumerate}
In particular, if $v$ is a valuation on $F$ such that $1+M \subseteq
\dot{F}^2$ and $B(F) = \dot{F}^2 U$, then an element $f \in F$ is double
rigid if and only if $v(f) \notin 2\Gamma$.
\end{Cor}

\begin{pf}
Since $F$ is nonrigid, there exists a valuation $v$ on $F$ such
that $1+M \subseteq \dot{F}^2$ and $B(F) = \dot{F}^2 U$.  We have $v$ is
not $2$-divisible if and only if $\dot{F}^2 U \subsetneq \dot{F}$,
which is equivalent to $B(F) \subsetneq \dot{F}$, in other words,
$F$ contains a double rigid element.
\end{pf}

We now define some of our main objects of study.

\begin{enumerate}
\item Let $F^{(2)} = F(\sqrt{F})$, the maximal multiquadratic extension of
$F$.  Then $\gal(F^{(2)}/F)$ is an elementary abelian $2$-group isomorphic
to $\prod_I {\Bbb Z}/2{\Bbb Z} $ where $|I| = \dim_{{\Bbb Z}/2{\Bbb Z}}
\dot{F}/\dot{F}^2$.  We denote $\gal(F^{(2)}/F)$ by $G_F^{[2]}$.

\item Let $F^{\{3\}} = F^{(2)}(\sqrt{F^{(2)}}) = (F^{(2)})^{(2)}$.  Then
$F^{\{3\}}$ is a Galois extension of $F$ and we denote $\gal(F^{\{3\}}/F)$
by $G_F^{\{3\}}$.

\item Let $F^{(3)}$ denote the subfield of $F^{\{3\}}$ generated by all
fields $F^{(2)}(\sqrt{g})$ such that $g \in F^{(2)}$ and
$F^{(2)}(\sqrt{g})/F$ is a Galois extension.  Denote $\gal(F^{(3)}/F)$ by
$G_F^{[3]}$.
\end{enumerate}

Thus $F \subseteq F^{(2)} \subseteq F^{(3)} \subseteq F^{\{3\}}$ and
$F^{(2)}/F$, $F^{(3)}/F$, $F^{\{3\}}/F$ are each Galois extensions.  The
groups
$G_F^{[3]}$, $G_F^{[2]}$ are quotients of $G_F^{\{3\}}$.

For the connection between $G_F^{[3]}$ and the Witt ring of quadratic forms
over $F$, see [MSp].

\begin{Lem}
\begin{enumerate}
\item The group $G_F^{\{3\}}$ has exponent dividing $4$.
\item For each $\sigma, \tau \in G_F^{\{3\}}$, the commutator $[\sigma,
      \tau]$ has order dividing $2$, and $[\sigma, \tau] = [\tau, \sigma]$.
\item All commutators and squares in $G_F^{\{3\}}$ commute with each other.
\end{enumerate}
\end{Lem}

\begin{pf}
Let $\sigma, \tau \in G_F^{\{3\}}$.  Then $\sigma^2|_{F^{(2)}} = 1$ since
$G_F^{[2]}$ has exponent dividing $2$ and $[\sigma, \tau]|_{F^{(2)}} =1$
since $G_F^{[2]}$ is abelian.  Thus $\sigma^2, [\sigma, \tau] \in
\gal(F^{\{3\}}/F^{(2)})$, a group of exponent dividing $2$.  Therefore
$\sigma^4 =1$ and $[\sigma, \tau]^2 =1$.  Note that $[\tau, \sigma] =
[\sigma, \tau]^{-1} = [\sigma, \tau]$  This proves (1) and (2).

For (3), we have from above that each commutator and square in
$G_F^{\{3\}}$ lies in $\gal(F^{\{3\}}/F^{(2)})$, which is an abelian group
of exponent dividing $2$.
\end{pf}

\begin{Lem}
For every $\sigma, \tau \in G_F^{\{3\}}$, the commutator
$[\, [\sigma, \tau], \,\sigma]$ commutes with both $\sigma$ and $\tau$.
\end{Lem}

\begin{pf}
In any group, we have the identity
\begin{displaymath}
[a,bc] =[a,c][a,b][\,[a,b],\,c].
\end{displaymath}
Let $b=c$ and assume that $a,b \in G_F^{\{3\}}$.  Then $[a,b^2]=
[\,[a,b],\,b]$ by Lemma 2.4.  Let $a =[\sigma, \tau]$ and $b=\sigma$.  We
have $[\,[\sigma, \tau], \,\sigma^2] =1$ by Lemma 2.4.  Therefore,
$[\,[\,[\sigma, \tau], \,\sigma], \,\sigma] =1$.  This implies
$[\,[\sigma, \tau], \,\sigma]$ commutes with $\sigma$.

Now let $a=\sigma^2$ and $b=\tau$.  Since $[\sigma^2, \tau^2]=1$ by Lemma
2.4, we have $[\,[\sigma^2, \tau],\,\tau]=1$ and so $[\sigma^2, \tau]$
commutes with $\tau$.  But
$[\sigma^2, \tau]= [\, [\sigma, \tau], \,\sigma]$ by the following argument.
Let $a =\tau$ and $b=\sigma$ in the identity to see
\begin{displaymath}
[\sigma^2, \tau]= [\tau, \sigma^2]= [\,[\tau, \sigma], \,\sigma] =
[\, [\sigma, \tau], \,\sigma].
\end{displaymath}
\end{pf}

\begin{remark}
The results in Lemmas 2.4, 2.5 hold in any homomorphic image of
$G_F^{\{3\}}$.
\end{remark}

\begin{Prop}  Every square and commutator of elements in $G_F^{[3]}$ lies in
$Z(G_F^{[3]})$, the center of $G_F^{[3]}$.
\end{Prop}

\begin{pf}
Let $\sigma, \tau, \psi \in G_F^{[3]}$.  We must show $\sigma^2$ and
$[\sigma, \tau]$ commute with $\psi$.  Since $F^{(3)}$ is the compositum of
all fields $F^{(2)}(\sqrt {g})$ where $g \in F^{(2)}$ and
$F^{(2)}(\sqrt {g})/F$ is a Galois extension, it is sufficient to prove the
result for elements of the group $\gal(F^{(2)}(\sqrt {g})/F)$ whenever
$F^{(2)}(\sqrt {g})/F$ is a Galois extension.  Since $G_F^{[2]}$ is an
abelian group of exponent dividing $2$, it follows $\sigma^2$ and $[\sigma,
\tau]$ both lie in $\gal(F^{(2)}(\sqrt{g})/F^{(2)})$, a group of order
dividing $2$.  As $F^{(2)}(\sqrt{g})/F$ is a Galois extension, the group
$\gal(F^{(2)}(\sqrt{g})/F^{(2)})$ is a normal subgroup of
$\gal(F^{(2)}(\sqrt{g})/F)$ of order $1$ or $2$ and thus lies in the center
of $\gal(F^{(2)}(\sqrt{g})/F)$.  Therefore $\sigma^2$ and $[\sigma, \tau]$
lie in $Z(\gal(F^{(2)}(\sqrt{g})/F))$.
\end{pf}

\begin{remark}
In Lemma 2.4(3) and Proposition 2.6, the result for commutators is a
consequence of the result for squares since the identity $[\sigma, \tau] =
(\sigma^{-1})^2 (\sigma \tau^{-1})^2 \tau^2$ shows that a commutator is
a product of three squares.
\end{remark}

\begin{Prop}
Let $G_0$ be a group with generators $x,y$ and assume
\begin{enumerate}
\item $x^4 =y^4=1$,
\item $[x,y]^2=[\,\,[x,y],x\,]^2 =[\,\,[x,y],y\,]^2 =1$,
\item $[\,\,[x,y],x\,]$ and $[\,\,[x,y],y\,]$ each commute with $x$ and
$y$.
\end{enumerate}
\noindent
Then the following identities hold in $G_0$, and $|G_0| \leq 128$.
\begin{enumerate}
\item[(a)] $yx =xy[x,y]$,
\item[(b)] $[x,y]x = x[x,y][\,\,[x,y],x\,]$,
\item[(c)] $[x,y]y = y[x,y][\,\,[x,y],y\,]$.
\end{enumerate}

\noindent
If in addition one assumes
\begin{enumerate}
\item[4.] $y^2 =1$,
\end{enumerate}
\noindent
then $[\,\,[x,y],y\,] =1$ and $|G_0| \leq 32$.
\end{Prop}

\begin{pf}
Since $[x,y]=[x,y]^{-1} = [y,x]$, we have $yx = xy[y,x]=xy[x,y]$, which is
(a).  For (b), let $z=[x,y]$ in the identity $zx=xz[z,x]$ and prove (c)
similarly.

We will next show that each element $g \in G_0$ can be written in the form
\[g= x^{e_1}y^{e_2}[x,y]^{e_3}[\,\,[x,y],x\,]^{e_4}[\,\,[x,y],y\,]^{e_5},\]
where $e_1,e_2 \in \{0,1,2,3\}$ and $e_3,e_4,e_5 \in \{0,1\}$.  This will
show $|G_0| \leq 4 \cdot 4 \cdot 2 \cdot 2 \cdot 2 =128$.  Let $g \in G_0$.
Since $x^4 =1$ and $[\,\,[x,y],x\,]$ commutes with $x$, (a) and (b) allow
us to write $g$ in the form $g=x^{e_1}z$ where $e_1 \in \{0,1,2,3\}$ and
$z$ involves factors of the type $y, [x,y], [\,\,[x,y],x\,]$.  Since
$y^4=1$ and $[\,\,[x,y],x\,]$ and $[\,\,[x,y],y\,]$ both commute with $y$,
(c) allows us to write $g=x^{e_1}y^{e_2}w$ where $e_2 \in \{0,1,2,3\}$ and
$w$ involves factors of the type $[x,y], [\,\,[x,y],x\,],
[\,\,[x,y],y\,]$.  We now use (2) and (3) to show $g$ has the required
form.

Now assume (4) also holds.  Then \[ [x,y]y= x^{-1}y^{-1}xy^2 = x^{-1}yx
=x^{-1}(xy[x,y]) =y[x,y] \] and so $[\,\,[x,y],y\,]=1$.  In this case we
may assume $e_2 \in \{0,1\}$ and $e_5 =0$, and therefore $|G_0| \leq 4 \cdot
2 \cdot 2 \cdot 2 \cdot 1 =32$.
\end{pf}

\begin{remark}
The proof that $[\,\,[x,y],y\,]=1$ in (4) holds in any group where $y^2 =
[x,y]^2 =1$.
\end{remark}

\begin{Cor}
\begin{enumerate}
\item Any subgroup of $G_F^{\{3\}}$ or any subgroup of a homomorphic image
of $G_F^{\{3\}}$ that is generated by two elements $x,y$ satisfies (1)-(3)
of Proposition 2.7.
\item Any subgroup of $G_F^{[3]}$ or any subgroup of a homomorphic image
of $G_F^{[3]}$ that is generated by two elements $x,y$ has order at most
$32$.
\end{enumerate}
\end{Cor}

\begin{pf}
(1).  Lemmas 2.4, 2.5 imply (1)-(3) hold in $G_F^{\{3\}}$, and these
properties are preserved in any homomorphic image of $G_F^{\{3\}}$.

(2).  As $G_F^{[3]}$ is a homomorphic image of $G_F^{\{3\}}$, we know
(1)-(3) hold.  Since $[x,y] \in Z(G_F^{[3]})$ by Proposition 2.6, it follows
$[\,\,[x,y],x\,] =[\,\,[x,y],y\,] =1$.  Therefore the subgroup has order at
most $4 \cdot 4 \cdot 2 \cdot 1 \cdot 1 = 32$.
\end{pf}

The next proposition and corollary are needed in section 4.

\begin{Prop}
Let $G$ be a metabelian group.  If $x,y,z \in G$ then
\[  [\,\,[x,y],\, z][\,\,[y,z],\, x][\,\,[z,x],\, y] =1. \]
\end{Prop}

\begin{pf}
Since $G$ is metabelian, there is a normal subgroup $H$ such that $H$ and
$G/H$ are abelian.  Thus every commutator of $G$ lies in $H$ and any two
commutators commute.

After completely expanding $[\,\,[x,y],\, z][\,\,[y,z],\,
x][\,\,[z,x],\, y]$ and making obvious cancellations, one is left with
\[ y^{-1}x^{-1}yxz^{-1}x^{-1}y^{-1}(xzyx^{-1}y^{-1}z^{-1})
(yxzy^{-1}z^{-1}x^{-1})zxy . \]
The two expressions in parentheses commute since they are the commutators
$[x^{-1}, y^{-1}z^{-1}]$ and $[y^{-1}, z^{-1}x^{-1}]$.  After transposing
these expressions and cancelling the remaining terms, one is left with the
identity.
\end{pf}

\begin{Cor}
Let $\sigma, \tau, \psi \in G_F^{\{3\}}$ and assume $[\tau, \psi] =1$.
Then $[\tau, \,\,[\sigma, \psi]\,]= [\psi, \,\,[\sigma, \tau]\,]$.
\end{Cor}

\begin{pf}
Since $G_F^{\{3\}}$ is metabelian, Proposition 2.9 implies
\[ 1= [\,\,[\sigma, \tau],\, \psi][\,\,[\tau, \psi],\, \sigma]
[\,\,[\psi, \sigma],\, \tau] =  [\,\,[\sigma, \tau],\, \psi] [\,\,[\psi,
\sigma],\, \tau]. \]
Since all commutators have order dividing $2$ by Lemma 2.4, the last
expression equals $[\psi, \,\,[\sigma, \tau]\,] [\tau, \,\,[\sigma,
\psi]\,]$ and the result follows.
\end{pf}

\section{The groups $D \curlywedge C$ and $G_1$ as Galois groups}

Let $D$ denote the dihedral group of order $8$ and let $C$ denote the
cyclic group of order $4$.  Let $\lambda : D \to {\Bbb Z}/2{\Bbb Z}$ be a
group homomorphism with $\ker(\lambda) \cong {\Bbb Z}/2{\Bbb Z} \times
{\Bbb Z}/2{\Bbb Z}$ and let $\eta : C \to {\Bbb Z}/2{\Bbb Z}$ be the unique
nontrivial group homomorphism.  We let $D \curlywedge C$ denote the
pullback of this pair of homomorphisms.  Thus $D \curlywedge C$ is the
subgroup of $D \times C$ of order $16$ consisting of elements $(u,v)$ such
that $\lambda(u) = \eta(v)$.  Note that the subgroup $D \curlywedge C$
projects onto both $D$ and $C$ with maps $\pi_1 : D \curlywedge C \to D$
and $\pi_2 : D \curlywedge C \to C$ such that $\lambda \circ \pi_1 = \eta
\circ \pi_2$.

Let $a,b \in F$ and let $E= F(\sqrt{a}, \sqrt{b})$.  Assume $[E:F]
=4$.  We let $D^{a,b}$ denote a Galois extension of $F$ (if one exists)
such that $F \subseteq E \subseteq D^{a,b}$, $\gal(D^{a,b}/F) \cong D$, and
$\gal(D^{a,b}/F(\sqrt{ab})) \cong C$.

We let $C^a$ denote a cyclic quartic extension of $F$ (if one exists) such
that $F \subseteq F(\sqrt{a}) \subseteq C^a$.

It is known that $D^{a,b}$ exists if and only if the quaternion algebra
$(a,b)_F = 0$ in the Brauer group $Br(F)$, and this is equivalent to $a \in
\im(N_{F(\sqrt{b})/F})$ where $N_{F(\sqrt{b})/F}$ is the norm map from
$F(\sqrt{b})$ to $F$.  The field $C^a$ exists if and only if the
quaternion algebra $(a,a)_F = 0$, which is equivalent to $a \in
\im(N_{F(\sqrt{a})/F})$ and also equivalent to $a$ being a sum of two
squares in $F$.

Now assume $C^a$ and $D^{a,b}$ exist, and let $K$ denote a composite of
$C^a$ and $D^{a,b}$.  We have $F(\sqrt{a}) = C^a \cap D^{a,b}$ since $D$
does not admit $C$ as a quotient group.  Thus $K/F$ is a Galois extension
and $[K:F]=16$.

\begin{Prop} $\gal(K/F) \cong D \curlywedge C$.
\end{Prop}

\begin{pf}
Let $r: \gal(K/F) \longrightarrow D \times C$ be the group
homomorphism defined by $r(\sigma) = (\sigma|_{D^{a,b}} , \sigma|_{C^a} )$.
Then $r$ is injective since $K= D^{a,b}C^a$.  Since $F(\sqrt{a}) = C^a \cap
D^{a,b}$ and $\gal(D^{a,b}/F(\sqrt{a})) \cong {\Bbb Z}/2{\Bbb Z} \times
{\Bbb Z}/2{\Bbb Z}$, it follows that $r(\sigma) \in D \curlywedge C$.
Since $|\gal(K/F)| = 16 = |D \curlywedge C|$, it follows that $r$ maps
$\gal(K/F)$ isomorphically onto $D \curlywedge C$.
\end{pf}

\begin{Prop}
The following statements are equivalent.
\begin{enumerate}
\item There is a Galois extension $K/F$ such that $\gal(K/F) \cong D
\curlywedge C$.
\item There exist elements $a,b \in F$ such that $[F(\sqrt{a}, \sqrt{b}) :
F] =4$ and $(a,a)_F = (a,b)_F =0 $ in the Brauer group $Br(F)$.
\end{enumerate}
\end{Prop}

\begin{pf}
We have that (2) implies (1) by Proposition 3.1 since (2) implies that
$D^{a,b}$ and $C^a$ exist.  Now suppose (1) holds and let $\gal(K/F) \cong
D \curlywedge C$.  The kernel of $\lambda \circ \pi_1 = \eta \circ \pi_2$
has order $8$ and thus the subfield of $K$ corresponding to $\ker(\lambda
\circ \pi_1)$ is a quadratic exension of $F$ that we denote $F(\sqrt{a})$.
It follows that $\ker(\pi_1)$ is a subgroup of $D \curlywedge C$ of order
$2$ which corresponds to a subfield of $K$ of the form $D^{a,b}$ for some
$b \in F$ and $\ker(\pi_2)$ is a subgroup of $D \curlywedge C$ of order $4$
which corresponds to a subfield of $K$ of the form $C^a$.  The existence of
both $D^{a,b}$ and $C^a$ implies (2).
\end{pf}

See [GSS] and [GS] for more information on the realizability of $D
\curlywedge C$ as a Galois group.

\begin{Cor} Let $K/F$ be a Galois extension with $\gal(K/F) \cong D
\curlywedge C$.  Then $K$ contains a unique quadratic extension of $F$ that
imbeds into a cyclic quartic extension contained in $K/F$.  This quadratic
extension imbeds into two different cyclic quartic extensions contained in
$K/F$.
\end{Cor}

\begin{pf}  This follows from the observation that $D \curlywedge C$ has
just two normal subgroups $N_1$, $N_2$ with the property that $(D
\curlywedge C)/N_i$ is a cyclic group of order $4$ and that $|N_1 \cap
N_2|=2$. Indeed if $D \curlywedge C /N \cong C$, then $N$ contains the
commutator subgroup $(D \curlywedge C)^{\prime}$.  Since $(D \times
C)^{\prime}$ has order $2$, it follows that $(D \curlywedge C)^{\prime}$
has order $2$ and $D \curlywedge C/(D \curlywedge C)^{\prime} \cong C
\times {\Bbb Z}/2{\Bbb Z}$.  One checks that $C \times {\Bbb Z}/2{\Bbb Z}$ has
exactly two subgroups $M$ of order $2$ such that
$(C \times {\Bbb Z}/2{\Bbb Z})/M \cong C$.  Then $N$ is one of the two
inverse images of $M$ in $D \curlywedge C$.  Since $|N_1 \cap N_2|= |(D
\curlywedge C)^{\prime}| = 2$, the subfields corresponding to $N_i$ do not
generate $K$, and so they must intersect in a quadratic extension of $F$.
\end{pf}

Now we consider a group $G_1$ (defined below) of order $32$ and
study how to contruct a Galois extension with Galois group isomorphic to
$G_1$.

\begin{Def}
Let $G_1$ be the group generated by two symbols $x,y$ subject to the
relations
\begin{enumerate}
\item $x^4 = y^2 = 1$, $[x,y]^2 =1$, $[\,\, [x,y],x \, ]^2=1$.
\item $[ [x,y],x ]$ commutes with $x$ and $y$.
\end{enumerate}
\end{Def}

\begin{Prop}
$|G_1| \leq 32$.
\end{Prop}

\begin{pf}
Since $[\,\, [x,y],y \, ]=1$ by Remark 3 after Proposition 2.7, the
result now follows from Proposition 2.7.
\end{pf}

We now give a construction that will show $|G_1| = 32$.
Let $a,b \in F$, let $E=F(\sqrt{a}, \sqrt{b})$ and assume $[E:F]=4$.
Assume in addition that $D^{a,b}$ and $C^a$ exist.  Then $a \in
\im(N_{F(\sqrt{b})/F}) \cap \im(N_{F(\sqrt{a})/F})$ and so there exist
$\alpha \in F(\sqrt{a})$ and $\beta \in F(\sqrt{b})$ such that
$N_{F(\sqrt{a})/F}(\alpha) = N_{F(\sqrt{b})/F}(\beta) =a$.  Then [Wd],
Lemma 2.14, implies there exist $\gamma \in E$ and $d \in F$ such that
$N_{E/F(\sqrt{a})}(\gamma) = \alpha d$ and $N_{E/F(\sqrt{b})}(\gamma) =
\beta d$.  Let $L = E(\sqrt{\alpha d}, \sqrt{\beta d}, \sqrt{\gamma})$ and
$K = E(\sqrt{\alpha d}, \sqrt{\beta d})$.

\begin{Prop}
\begin{enumerate}
\item $K/F$ is a Galois extension with $\gal(K/F) \cong D \curlywedge C$.
\item $L/F$ is a Galois extension with $\gal(L/F) \cong G_1$.
\item $K \subseteq F^{(3)}$, $L \subseteq F^{\{3\}}$, $L \nsubseteq
F^{(3)}$.
\end{enumerate}
\end{Prop}

\begin{pf}
The extensions $F(\sqrt{\alpha d})/F$ and $E(\sqrt{\beta d}) /F$ are Galois
with
\[\gal(F(\sqrt{\alpha d})/F) \cong C, \,\,\,\, \gal(E(\sqrt{\beta d}) /F)
\cong D,\]
because $N_{F(\sqrt{a})/F}(\alpha d) = ad^2$ and
$N_{F(\sqrt{b})/F}(\beta d) = ad^2$.  We also note that $\gal(E(\sqrt{\beta
d})/F(\sqrt{ab})) \cong C$ because $F(\sqrt{ab}) \subseteq
F(\sqrt{ab})(\sqrt{a})=E \subseteq E(\sqrt{\beta d})$ and
$N_{E/F(\sqrt{ab})}(\beta d) = N_{F(\sqrt{b})/F}(\beta d) = ad^2$.  We see
that $K$ is the composite of $F(\sqrt{\alpha d})$ and $E(\sqrt{\beta
d})$, so $K/F$ is Galois and Proposition 3.1 implies $\gal(K/F) \cong D
\curlywedge C$.

We have $K \subseteq F^{(3)}$ since $F(\sqrt{\alpha d})=
F(\sqrt{a})(\sqrt{\alpha d}) \subseteq F^{(3)}$ and $E(\sqrt{\beta d})
\subseteq F^{(3)}$.  This proves (1) and the first part of (3).

Since $L=K(\sqrt{\gamma})$, in order to show $L/F$ is Galois, it is
sufficient to show $\gamma \sigma(\gamma) \in K^2$ for all $\sigma \in
\gal(K/F)$.  Since $\gamma \in E$, we need only consider $\sigma \in
\gal(E/F)$.  This is clear for $\sigma = 1$, and the others follow from the
calculations $N_{E/F(\sqrt{a})}(\gamma) = \alpha d \in K^2$,
$N_{E/F(\sqrt{b})}(\gamma) = \beta d \in K^2$ and
\begin{displaymath}
N_{E/F(\sqrt{ab})}(\gamma) = \frac{N_{E/F}(\gamma)
\gamma^2}{N_{E/F(\sqrt{a})}(\gamma) N_{E/F(\sqrt{b})}(\gamma)} = \frac{ad^2
\gamma^2}{\alpha d \beta d} \in K^2.
\end{displaymath}
Therefore $L/F$ is a Galois extension.

We have $L\subseteq F^{\{3\}}$ since $\alpha d$, $\beta d$, $\gamma$ $\in E
\subseteq F^{(2)}$.  To show $L \nsubseteq F^{(3)}$ it is sufficient to
show $F^{(2)}(\sqrt{\gamma})/F$ is not a Galois extension.  There exists
$\sigma \in \gal(F^{(2)}/F)$ such that $\gamma \sigma(\gamma) =
N_{E/F(\sqrt{a})}(\gamma) = \alpha d$.  But $\alpha d \notin (F^{(2)})^2$
since $\gal(F(\sqrt{\alpha d})/F) \cong C$, and thus
$F^{(2)}(\sqrt{\gamma})/F$ is not a Galois extension.  This finishes the
proof of (3).

We now show $\gal(L/F) \cong G_1$.  We have $[L:F]=32$ since $[K:F]=16$ and
$K \subsetneq L$ by (3).  There exist $\sigma_a, \sigma_b \in \gal(L/F)$
such that $\sigma_a(\sqrt{a}) =-\sqrt{a}$, $\sigma_a(\sqrt{b})= \sqrt{b}$,
$\sigma_b(\sqrt{a}) =\sqrt{a}$, $\sigma_b(\sqrt{b}) =-\sqrt{b}$ and
$\sigma_b(\sqrt{\alpha d}) = \sqrt{\alpha d}$.  (To see that $\sigma_b$
exists, note that an automorphism with this property lies in $\gal(K/F)$
and that it can be extended to $\gal(L/F)$.)  The group generated by
$\sigma_a$ and $\sigma_b$ equals $\gal(L/F)$, otherwise the fixed field of this
subgroup would contain a quadratic extension of $F$, but none of $\sqrt{a},
\sqrt{b}, \sqrt{ab}$ is fixed by both $\sigma_a$ and $\sigma_b$.
We will now show that $\sigma_a$ and $\sigma_b$ satisfy the relations
in Definition 3.4.  Since the group generated by $\sigma_a$ and $\sigma_b$
has order $32$ and $|G_1| \leq 32$, this will imply that $\gal(L/F) \cong
G_1$.

Lemmas 2.4, 2.5 and Remark 1 following these lemmas, along with the fact
that $\sigma_a$ and $\sigma_b$ generate $\gal(L/F)$, imply that $\sigma_a^4
= [\sigma_a, \sigma_b]^2 = [ [\sigma_a, \sigma_b], \sigma_a]^2 = 1$, and $[
[\sigma_a, \sigma_b], \sigma_a] \in Z(G_1)$.  It remains to prove
$\sigma_b^2 =1$.

Since we know $\gal(E(\sqrt{\beta d})/F) \cong D$, $\gal(E(\sqrt{\beta
d})/F(\sqrt{ab})) \cong C$ and $\sigma_b(\sqrt{ab}) = - \sqrt{ab}$, it
follows that $\sigma_b|_{E(\sqrt{\beta d})}$ has order $2$.
Thus $\sigma_b|_K$ has order $2$.  We now show $\sigma_b^2(\sqrt{\gamma}) =
\sqrt{\gamma}$.
We have $\gamma \sigma_b(\gamma) = \alpha d$ and thus
$\sqrt{\gamma} \sigma_b(\sqrt{\gamma}) = (-1)^{\epsilon} \sqrt{\alpha d}$,
where $\epsilon \in \{0,1\}$.  Then
$\sigma_b(\sqrt{\gamma}) \sigma_b^2(\sqrt{\gamma}) =
\sigma_b((-1)^{\epsilon} \sqrt{\alpha d}) =
(-1)^{\epsilon} \sqrt{\alpha d}$.
These equations imply $\sigma_b^2(\sqrt{\gamma}) = \sqrt{\gamma}$.

It follows that $\gal(L/F) \cong G_1$ since the relations in Definition 3.4
hold and $[L:F]=32$.
\end{pf}

\begin{Thm}  The following statements are equivalent.
\begin{enumerate}
\item There is a Galois extension $K/F$ such that $\gal(K/F) \cong D
\curlywedge C$.
\item There is a Galois extension $L/F$ such that $\gal(L/F) \cong G_1$.
\item There exist $a,b \in F$ such that $[F(\sqrt{a}, \sqrt{b}):F]=4$ and
$(a,a)_F = (a,b)_F =0$.
\item There exists $a \in D_F(\langle 1,1 \rangle) \backslash \dot{F}^2$
such that $|D_F(\langle 1, -a \rangle)| \geq 4$.
\item $F$ is not rigid and $D_F(\langle 1,1 \rangle) \nsubseteq F^2 \cup
-F^2$.
\item Either $-1 \in F^2$ and $F$ is not rigid or $-1 \notin F^2$ and
$D_F(\langle 1,1 \rangle) \nsubseteq F^2 \cup -F^2$.
\end{enumerate}

If $F$ is formally real, then statements (1)-(6) hold if and only if $F$ is
not a Pythagorean field.
\end{Thm}

\begin{pf}
We have already seen in Proposition 3.2 and the proof of Proposition 3.6 that
(1)-(3) are equivalent.  Note that (2) implies (1) because $D \curlywedge
C$ is a quotient of $G_1$.  The equivalence of (3) and (4) is easy to
check.  It is clear that (5) implies (6).

(6) $\Rightarrow$ (3):  First assume $-1 \in F^2$ and $F$ is not rigid.
Then there exists $a \in F$ with $a \notin \pm F^2$ such that $b \in
D_F(\langle 1,a \rangle)$ and $b \notin F^2 \cup aF^2$.  Then (3) holds
since $(a,b)_F =0$ and $-1 \in F^2$ implies $(a,a)_F = (a, -a)_F =0$.  Now
assume $-1 \notin F^2$ and $D_F(\langle 1,1 \rangle) \nsubseteq F^2 \cup
-F^2$.  Let $a \in D_F(\langle 1,1 \rangle)$ with $a \notin \pm F^2$.  Then
$[F(\sqrt{a}, \sqrt{-1}):F] =4$ and $(a,a)_F =0$.  Now (3) holds with $b =
-a$.

(3) $\Rightarrow$ (5): Suppose $D_F(\langle 1,1 \rangle) \subseteq F^2 \cup
-F^2$.  We have $a \in D_F(\langle 1,1 \rangle)$ since $(a,a)_F =0$.  Since
$a \notin F^2$, we have $a \in -F^2$.  Then $(-1,b)_F =0$ and thus $b \in
D_F(\langle 1,1 \rangle)$.  This is a contradiction since $b \notin F^2
\cup aF^2 = F^2 \cup -F^2$.  Therefore $D_F(\langle 1,1 \rangle) \nsubseteq
F^2 \cup -F^2$.  Since $(a,a)_F = (a,b)_F =0$ implies $a,b \in D_F(\langle
1, -a \rangle)$, we have $1,a,b,ab \in D_F(\langle 1, -a \rangle)$ and so
$D_F(\langle 1, -a \rangle)$ contains at least four square classes. Thus
$F$ is not rigid by Lemma 2.1.

Now assume $F$ is formally real.  If (6) holds, then $-1 \notin F^2$ and
$D_F(\langle 1,1 \rangle) \nsubseteq F^2 \cup -F^2$.  Thus $F$ is not a
Pythagorean field.  If $F$ is not a Pythagorean field, then $D_F(\langle
1,1 \rangle) \nsubseteq F^2 \cup -F^2$, since $F$ is formally real.  Thus
(6) holds.
\end{pf}

Nondyadic local fields are rigid.  Every dyadic local field $F$ satisfies
conditions (1)-(6).  In fact, $a$ can be chosen as any nonsquare element in
$D_F(\langle 1,1 \rangle)$.  All global fields (of characteristic different
from $2$) also satisfy conditions (1)-(6).

\begin{Prop}
Suppose $\gal(L/F) \cong G_1$.  Then there is a unique quadratic extension
$F(\sqrt{a})$ of $F$ such that
\[ F \subseteq F(\sqrt{a}) \subseteq E \subseteq L, \]
for some subfield $E$ of $L$ where $\gal(E/F) \cong C$.  The element $-a$
is not rigid.
\end{Prop}

\begin{pf}
Let $K = F^{(3)} \cap L$.  Then $\gal(K/F) \cong D \curlywedge C$.
Corollary 3.3 implies there is a unique quadratic extension $F(\sqrt{a})$
of $F$ in $K$ such that $F$ is contained in a cyclic quartic extension $E$
lying in $K$.  The construction in Proposition 3.1 along with
Proposition 3.2 and Theorem 3.7 imply $-a$ is not rigid.
\end{pf}

In the notation of this section, the two cyclic quartic extensions that
appear in Corollary 3.3 are easily seen to be $F(\sqrt{a}, \sqrt{\alpha
d})$ and $F(\sqrt{a}, \sqrt{\alpha db})$.

\section{The groups $D \curlywedge D$ and $G_2$ as Galois groups}

Let $\lambda_1 : D \to {\Bbb Z}/2{\Bbb Z}$ be a group homomorphism with
$\ker(\lambda_1) \cong {\Bbb Z}/2{\Bbb Z} \times {\Bbb Z}/2{\Bbb Z}$ and
let $\lambda_2 : D \to {\Bbb Z}/2{\Bbb Z}$ be another group homomorphism
with $\ker(\lambda_2) \cong {\Bbb Z}/2{\Bbb Z} \times {\Bbb Z}/2{\Bbb Z}$.
Let $D \curlywedge D$ denote the pullback of this pair of homomorphisms.
Thus $D \curlywedge D$ is the subgroup of $D \times D$ of order $32$
consisting of elements $(u,v)$ such that $\lambda_1(u) = \lambda_2(v)$.
Note that the subgroup $D \curlywedge D$ projects onto $D$ with maps
$\pi_1 : D \curlywedge D \to D$ and $\pi_2 : D \curlywedge D \to D$ such
that $\lambda_1 \circ \pi_1 = \lambda_2 \circ \pi_2$.

\begin{Prop}
The following statements are equivalent.
\begin{enumerate}
\item There is a Galois extension $K/F$ such that $\gal(K/F) \cong D
\curlywedge D$.
\item There exist elements $a,b,c \in F$ such that $[F(\sqrt{a}, \sqrt{b},
\sqrt{c}): F] =8$ and $(a,b)_F = (a,c)_F =0 $ in the Brauer group $Br(F)$.
\end{enumerate}
\end{Prop}

\begin{pf}
First assume that (2) holds.  Then there exist Galois extensions $D^{a,b}$
and $D^{a,c}$ of $F$.  Let $K$ denote a composite of $D^{a,b}$
and $D^{a,c}$.  Then $D^{a,b} \cap D^{a,c} = F(\sqrt{a})$ since the
intersection is a Galois extension of $F$ and the only Galois subextension
of $D^{a,b}$ of degree $4$ over $F$ is $F(\sqrt{a}, \sqrt{b})$.  Therefore,
$K/F$ is a Galois extension of $F$ with $[K:F]=32$.  Now a proof that is
very similar to the proof of Proposition 3.1 shows that $\gal(K/F) \cong D
\curlywedge D$.

Now assume that (1) holds and let $\gal(K/F) \cong D \curlywedge D$.  The
kernel of $\lambda_1 \circ \pi_1 = \lambda_2 \circ \pi_2$ has order $16$
and thus the subfield of $K$ corresponding to this kernel is a quadratic
extension of $F$ which we denote $F(\sqrt{a})$.  Since both $\pi_1$ and
$\pi_2$ are surjective, it follows that $\ker(\pi_1)$ and $\ker(\pi_2)$ are
each normal subgroups of $D \curlywedge D$ of order $4$ which correspond to
subfields $E_1, E_2$ of $K$ which are Galois over $F$ with $\gal(E_1/F)
\cong D \cong \gal(E_2/F)$.  Since $\ker(\pi_1) \cap \ker(\pi_2) =1$, it
follows that $\ker(\pi_1)$ and $\ker(\pi_2)$ generate a subgroup of order
$16$.  This subgroup is $\ker(\lambda_1 \circ \pi_1)$ since $\ker(\pi_i)
\subseteq \ker(\lambda_1 \circ \pi_1) = \ker(\lambda_2 \circ \pi_2)$.  Thus
$E_1 \cap E_2 = F(\sqrt{a})$.  Since $\ker(\lambda_i) \cong {\Bbb Z}/2{\Bbb
Z} \times {\Bbb Z}/2{\Bbb Z}$, it follows $E_1$ has the form $D^{a,b}$ and
$E_2$ has the form $D^{a,c}$ for some $b,c \in \dot{F}$.  Thus $(a,b)_F =
(a,c)_F =0$.  We have $[F(\sqrt{a}, \sqrt{b}, \sqrt{c}): F] =8$ since
$D^{a,b} \cap D^{a,c} = F(\sqrt{a})$.
\end{pf}

\begin{Prop}
The group $D \curlywedge D$ contains a unique abelian subgroup of order
$16$.  This subgroup equals $\ker(\lambda_1) \times \ker(\lambda_2)$ and is
isomorphic to $({\Bbb Z}/2{\Bbb Z})^4$.
\end{Prop}

\begin{pf}
We have that $\ker(\lambda_1) \times \ker(\lambda_2)$ is an abelian
subgroup of $D \curlywedge D$ of order $16$ that is isomorphic to $({\Bbb
Z}/2{\Bbb Z})^4$.  Now let $H$ be any abelian subgroup of $D \curlywedge D$
of order $16$.  Since $\pi_i(H)$ is an abelian subgroup of $D$, we have
$\pi_i(H)$ is a subgroup of order at most $4$.  Since $H \subseteq \pi_1(H)
\times \pi_2(H)$ and $H$ has order $16$, it follows $H = \pi_1(H) \times
\pi_2(H)$.  Since $\pi_1(H) \times \pi_2(H) \subseteq D \curlywedge D$, it
follows $\pi_1(H) = \ker(\lambda_1)$ and $\pi_2(H) = \ker(\lambda_2)$.
\end{pf}

\begin{Cor}
In the notation of Proposition 4.1,
\[\gal(K/F(\sqrt{a})) = \ker(\lambda_1 \circ \pi_1) = \ker(\lambda_1)
\times \ker(\lambda_2) \cong ({\Bbb Z}/2{\Bbb Z})^4. \]
\end{Cor}

\begin{pf}
The first equality follows from the proof of Proposition 4.1.  Since
$\ker(\lambda_1) \times \ker(\lambda_2) \subseteq \ker(\lambda_1 \circ
\pi_1)$ and each has order $16$, we have equality.
\end{pf}

Now we consider a group $G_2$ (defined below) of order $64$ and
study how to construct a Galois extension over a field $F$ with Galois
group isomorphic to $G_2$.

\begin{Def}
Let $G_2$ be the group generated by three symbols $x,y,z$ subject to the
relations
\begin{enumerate}
\item $x^2=y^2=z^2=1$, $[x,y]^2 =[x,z]^2 = [y,z] =1$.
\item $[y, [x,z] ]= [z, [x,y] ]$ has order dividing 2 and commutes with
$x$, $y$ and $z$.
\end{enumerate}
\end{Def}

\begin{Lem}
The following identities hold in $G_2$ and $|G_2| \leq 64$.
\begin{enumerate}
\item $yx =xy[x,y]$, $zx=xz[x,z]$.
\item $[x,y]$ commutes with both $x$ and $y$ and $[x,z]$ commutes with both
$x$ and $z$.
\item $[x,z]y =y[x,z][y, \,[x,z]\,]$ and $[x,y]z =z[x,y][z, \,[x,y]\,]$.
\item $[x,z][x,y] = [x,y][x,z]$.
\end{enumerate}
\end{Lem}

\begin{pf}
Statements (1)-(3) are proved as in Proposition 2.7.  (See Remark 3
following the proof of Proposition 2.7.)  Since $[x,y] = xyxy$, and
$[y, [x,z] ]$ is in the center, we have
\begin{align*}
\begin{split}
[x,z][x,y] &= [x,z]xyxy =x[x,z]yxy =xy[x,z][y,\, [x,z]\, ]xy \\
 &= xy[x,z]xy[y,\, [x,z]\, ] = xyxy[x,z][y,\, [x,z]\, ]^2 \\
 &= [x,y][x,z],
\end{split}
\end{align*}
which proves (4).

To show $|G_2| \leq 64$, it is sufficient to show that each element of
$G_2$ can be written in the form
\begin{displaymath}
x^{e_1}y^{e_2}z^{e_3}[x,y]^{e_4}[x,z]^{e_5}[y,\, [x,z]\, ]^{e_6},
\end{displaymath}
where each $e_i \in \{0,1\}$, $1 \leq i \leq 6$.  The relations defining
$G_2$ and statements (1)-(4) allow exactly this.
\end{pf}

We now give a construction of a Galois extension $L/F$ with $[L:F] = 64$
and $\gal(L/F) \cong G_2$.

Assume $a,b,c \in F$ such that $[F(\sqrt{a}, \sqrt{b}, \sqrt{c}):F] = 8$
and $(a,b)_F = (a,c)_F =0$ in the Brauer group $Br(F)$.  Then $a \in
\im(N_{F(\sqrt{b})/F}) \cap \im(N_{F(\sqrt{c})/F})$ and so there exist
$\beta \in F(\sqrt{b})$ and $\gamma \in F(\sqrt{c})$ such that
$N_{F(\sqrt{b})/F}(\beta) = N_{F(\sqrt{c})/F}(\gamma) =a$.
Let $E=F(\sqrt{b}, \sqrt{c})$.  Then [Wd], Lemma 2.14, implies there exist
$\delta \in E$ and $d \in F$ such that $N_{E/F(\sqrt{b})}(\delta) = \beta
d$ and $N_{E/F(\sqrt{c})}(\delta) = \gamma d$.

Let $E^{\prime} = E(\sqrt{a})$, $K=E^{\prime}(\sqrt{\beta d}, \sqrt{\gamma
d})$, and $L=K(\sqrt{\delta})$.

\begin{Prop}
\begin{enumerate}
\item $K/F$ is a Galois extension with $\gal(K/F) \cong D \curlywedge D$.
\item $L/F$ is a Galois extension with $\gal(L/F) \cong G_2$.
\item $K \subseteq F^{(3)}$, $L \subseteq F^{\{3\}}$, $L \nsubseteq
F^{(3)}$.
\end{enumerate}
\end{Prop}

\begin{pf}
We have $F(\sqrt{a}, \sqrt{b}, \sqrt{\beta d}) = D^{a,b}$ and $F(\sqrt{a},
\sqrt{c}, \sqrt{\gamma d}) = D^{a,c}$.  Since $D^{a,b}$ and $D^{a,c}$ are
each Galois extensions of $F$, it follows that $K= D^{a,b}D^{a,c}$ is a
Galois extension of $F$.  The reasoning in the proof of Proposition 4.1
lets us conclude $\gal(K/F) \cong D \curlywedge D$.  This proves (1).

We have $D^{a,b}$ and $D^{a,c}$ contained in $F^{(3)}$ and thus $K
\subseteq F^{(3)}$.  Since $L= E^{\prime}(\sqrt{\beta d}, \sqrt{\gamma d},
\sqrt{\delta})$ and $\beta d$, $\gamma d$, $\delta$ $\in F^{(2)}$, it
follows $L \subseteq F^{\{3\}}$.  If $L \subseteq F^{(3)}$, then
$F^{(2)}(\sqrt{\delta})/F$ is a Galois extension.  This implies
$N_{E/F(\sqrt{b})}(\delta) = \beta d \in (F^{(2)})^2$.  But this is
impossible since $D^{a,b} \nsubseteq F^{(2)}$.  This proves (3).

Since $[L:K]=2$ by (3) and $[K:F]=32$ by (1), it follows $[L:F]=64$.
Since $\delta \in E$ and $K/F$ is a Galois extension, in order to show
$L/F$ is a Galois extension, it is sufficient to show $\sigma(\delta)\delta
\in K^2$ for all $\sigma \in \gal(E/F)$.  This is accomplished as in the
proof of Proposition 3.6.

Now we begin to show $\gal(L/F) \cong G_2$.
There exist automorphisms $\sigma_a, \sigma_b, \sigma_c \in \gal(L/F)$ such
that $\sigma_a$ fixes $\sqrt{b}, \sqrt{c}$ and $\sigma_a(\sqrt{a}) =
-\sqrt{a}$, and similarly $\sigma_b$ fixes $\sqrt{a}, \sqrt{c}$, but not
$\sqrt{b}$, and $\sigma_c$ fixes $\sqrt{a}, \sqrt{b}$, but not
$\sqrt{c}$.

We want to choose $\sigma_b, \sigma_c $ more carefully as follows.  Choose
$\sigma_b \in \gal(L/F)$  such that
\[ \sigma_b|_{F(\sqrt{a}, \sqrt{c}, \sqrt{\gamma d})} =1, \,\,\,\,
\sigma_b^2|_{F(\sqrt{a}, \sqrt{b}, \sqrt{\beta d})} =1, \,\,\,\,
\sqrt{\beta d} \,\sigma_b(\sqrt{\beta d}) = \sqrt{a} d. \]
To see this is possible, observe that $\gal(F(\sqrt{a}, \sqrt{b},
\sqrt{\beta d})/F) \cong D$ and hence there exists $\sigma_b$ such that the
second and third conditions hold.  Since $K = F(\sqrt{a}, \sqrt{b},
\sqrt{\beta d}) F(\sqrt{a}, \sqrt{c}, \sqrt{\gamma d})$, $\gal(K/F) \cong D
\curlywedge D$ and $\sigma_b|_{F(\sqrt{a})} = 1$, the first condition can
be arranged simultaneously with the second and third.

Similarly we can choose $\sigma_c \in \gal(L/F)$ such that
\[ \sigma_c|_{F(\sqrt{a}, \sqrt{b}, \sqrt{\beta d})} =1, \,\,\,\,
\sigma_c^2|_{F(\sqrt{a}, \sqrt{c}, \sqrt{\gamma d})} =1, \,\,\,\,
\sqrt{\gamma d} \,\sigma_c(\sqrt{\gamma d}) = \sqrt{a} d. \]

We now show $\sigma_a, \sigma_b, \sigma_c $ satisfy the relations in
Definition 4.4.
Since $\sigma_a$ fixes each element of $E$ and $\delta \in E$, it follows that
$\sigma_a^2(\sqrt{\delta}) = \sqrt{\delta}$.  Similar reasoning shows that
$\sigma_a^2$ fixes $E^{\prime}$ and $K$ elementwise and hence $\sigma_a^2
=1$.

We have that $\sigma_b^2|_K = \sigma_c^2|_K = 1$.  To show $\sigma_b^2 =
\sigma_c^2 =1$, it is enough to show $\sigma_b^2(\sqrt{\delta}) =
\sqrt{\delta}$ and $\sigma_c^2(\sqrt{\delta}) = \sqrt{\delta}$.
Since
\[ \delta\,\sigma_b(\delta) = N_{E/F(\sqrt{c})}(\delta) =\gamma d, \]
we have $\sqrt{\delta} \,\sigma_b(\sqrt{\delta})=
(-1)^{\epsilon}\sqrt{\gamma d}$ and thus
\[\sigma_b(\sqrt{\delta})\, \sigma_b^2(\sqrt{\delta}) =
(-1)^{\epsilon}\sigma_b(\sqrt{\gamma d}) = (-1)^{\epsilon}\sqrt{\gamma d}. \]
This implies $\sigma_b^2(\sqrt{\delta}) =
\sqrt{\delta}$.  Since $\delta\,\sigma_c(\delta) =
N_{E/F(\sqrt{b})}(\delta) =\beta d$, we have $\sqrt{\delta}
\,\sigma_c(\sqrt{\delta})= (-1)^{\epsilon^{\prime}}\sqrt{\beta d}$, and a
similar proof gives $\sigma_c^2(\sqrt{\delta}) = \sqrt{\delta}$.

We have $[\sigma_a, \sigma_b]^2 = [\sigma_a, \sigma_c]^2 = 1$ since
$\gal(L/F)$ is a homomorphic image of $G_F^{\{3\}}$.  We have
$[\sigma_b, \sigma_c]|_K =1$ because $\sigma_b\sigma_c = \sigma_c\sigma_b$
on each subfield $F(\sqrt{a}, \sqrt{b}, \sqrt{\beta d})$ and
$F(\sqrt{a}, \sqrt{c}, \sqrt{\gamma d})$.  To show $[\sigma_b, \sigma_c]
=1$, it is enough to show $\sigma_b\sigma_c(\sqrt{\delta}) =
\sigma_c\sigma_b(\sqrt{\delta})$.

The equation $\sqrt{\delta} \,\sigma_b(\sqrt{\delta})=
(-1)^{\epsilon}\sqrt{\gamma d}$ gives
\[\sigma_c(\sqrt{\delta})\, \sigma_c\,\sigma_b(\sqrt{\delta}) =
(-1)^{\epsilon}\sigma_c(\sqrt{\gamma d}),\]
and the equation $\sqrt{\delta} \,\sigma_c(\sqrt{\delta})=
(-1)^{\epsilon^{\prime}}\sqrt{\beta d}$ gives
\[\sigma_b(\sqrt{\delta})\, \sigma_b\,\sigma_c(\sqrt{\delta}) =
(-1)^{\epsilon^{\prime}}\sigma_b(\sqrt{\beta d}).\]
Since $\sqrt{\beta d} \,\sigma_b(\sqrt{\beta d}) = \sqrt{\gamma d}
\,\sigma_c(\sqrt{\gamma d})$, some routine manipulation gives
$\sigma_b\sigma_c(\sqrt{\delta}) = \sigma_c\sigma_b(\sqrt{\delta})$.

We now show $[\sigma_b, \,\,[\sigma_a, \sigma_c]\,] =
[\sigma_c, \,\,[\sigma_a, \sigma_b]\,]$ has order dividing $2$ and commutes
with $\sigma_a, \sigma_b, \sigma_c$.  The two expressions are equal by
Corollary 2.10.  The order divides $2$ by Lemma 2.4.
We have $[\sigma_a, \sigma_b]$ in the center of $\gal(F(\sqrt{a}, \sqrt{b},
\sqrt{\beta d})/F)$ by Proposition 2.6 and thus $[\sigma_c, \,\,[\sigma_a,
\sigma_b]\,]| _{F(\sqrt{a}, \sqrt{b}, \sqrt{\beta d})}= 1$.  Similarly,
$[\sigma_a, \sigma_c]$ is in the center of $\gal(F(\sqrt{a}, \sqrt{c},
\sqrt{\gamma d})/F)$ and thus $[\sigma_b, \,\,[\sigma_a, \sigma_c]\,]|
_{F(\sqrt{a}, \sqrt{c}, \sqrt{\gamma d})}= 1$.  Thus, $[\sigma_b,
\,\,[\sigma_a, \sigma_c]\,]|_K =1$ and so
\[ [\sigma_b, \,\,[\sigma_a, \sigma_c]\,] \in \gal(L/K) \subseteq
Z(\gal(L/F)).\]
Therefore $[\sigma_b, \,\,[\sigma_a, \sigma_c]\,]$ commutes with $\sigma_a,
\sigma_b, \sigma_c$.

We have now shown that $\sigma_a, \sigma_b, \sigma_c$ satisfy the relations
in Definition 4.4.  Next we show $\sigma_a, \sigma_b, \sigma_c$ generate
$\gal(L/F)$.

Let $H = \langle \sigma_a, \sigma_b, \sigma_c \rangle$. If $H \subsetneq
\gal(L/F)$, then there exists a quadratic extension $F(\sqrt{e})$ of $F$ in
the fixed field of $H$ and $F(\sqrt{e}) \subseteq F^{(3)} \cap L =K$.  An
easy check shows that no quadratic extension of $F$ inside $F(\sqrt{a},
\sqrt{b}, \sqrt{c})$ is in the fixed field of $H$.  Then $F(\sqrt{a},
\sqrt{b}, \sqrt{c}, \sqrt{e}) \subseteq K$ and $({\Bbb Z}/2{\Bbb Z})^4$ is
a quotient of $\gal(K/F) \cong D \curlywedge D$ by a normal subgroup $N$ of
order $2$.  This implies the commutator subgroup of $D \curlywedge D$ is
contained in $N$ and this is impossible since the commutator subgroup of $D
\curlywedge D$ has order $4$.  Thus $({\Bbb Z}/2{\Bbb Z})^4$ is not a
quotient of $D \curlywedge D$ and therefore $H= \gal(L/F)$.  Since
$\sigma_a, \sigma_b, \sigma_c$ generate $\gal(L/F)$ and satisfy the
relations in Definition 4.4, it follows that $\gal(L/F) \cong G_2$.
\end{pf}

\begin{Thm}
The following statements are equivalent.
\begin{enumerate}
\item There is a Galois extension $K/F$ such that $\gal(K/F) \cong D
\curlywedge D$.
\item There is a Galois extension $L/F$ such that $\gal(L/F) \cong G_2$.
\item There exist $a,b,c \in F$ such that
$[F(\sqrt{a}, \sqrt{b}, \sqrt{c}):F]=8$ and
$(a,b)_F = (a,c)_F =0$.
\item There exist $a,b,c \in \dot{F}\backslash \dot{F}^2$ such that $a,b,c$
are independent mod squares and $b,c \in D_F(\langle 1, -a \rangle)$.
\item There exists $a \in \dot{F}\backslash \dot{F}^2$ such that either
$|D_F(\langle 1, -a \rangle)| \geq 8$ or both $|D_F(\langle 1, -a \rangle)|
=4$ and $a \notin D_F(\langle 1, 1  \rangle)$.
\end{enumerate}
\end{Thm}

\begin{pf}
(1) and (3) are equivalent by Proposition 4.1.  The construction in
Proposition 4.6 shows (3) implies (2).  (2) implies (1) since $D
\curlywedge D$ is a quotient of $G_2$ (as the construction in Proposition
4.6 shows).  It is straightforward to check that (3)-(5) are equivalent
using the observation $a \in D_F(\langle 1, -a \rangle)$ if and only if
$a \in D_F(\langle 1, 1  \rangle)$.
\end{pf}

In the following proposition, $DD$ denotes the central product of two copies
of the dihedral group $D$.  That is, if $g$ denotes the nontrivial element
in the center of $D$, then $DD \cong (D \times D)/\{1, (g,g)\}$.  Thus $DD$
is obtained from the direct product $D \times D$ by identifying the
centers of the two copies of $D$.

\begin{Prop}
Suppose $\gal(L/F) \cong G_2$.  Let $K = F^{(3)} \cap L$, so that
$\gal(K/F) \cong D \curlywedge D$.  Then there is a unique quadratic extension
$F(\sqrt{a})$ of $F$ such that $\gal(K/F(\sqrt{a}))$ is an
abelian group of order $16$.  In addition the following hold.
\begin{enumerate}
\item There is a unique subgroup $N$ of $G_2$ of order $2$ such that
$G_2/N \cong D \curlywedge D$.
\item $\gal(L/F(\sqrt{a})) \cong DD$.
\item $\gal(L/F(\sqrt{a}))$ is the unique subgroup of $G_2$ isomorphic to
$DD$.
\item The element $-a$ is nonrigid.
\end{enumerate}
\end{Prop}

\begin{pf}
Using Proposition 4.2, let $F(\sqrt{a})$ be the quadratic extension of $F$
that corresponds to the unique abelian subgroup of $\gal(K/F)$ of order
$16$.  Any extension $K'$ of $F$ that lies in $L$ with $\gal(K'/F)
\cong D \curlywedge D$ must lie in $F^{(3)}$.  Thus $K$ is uniquely
determined, so this proves (1).

In the notation of Proposition 4.6, let $D_1 = \langle \sigma_b, [\sigma_a,
\sigma_c] \rangle$, let $D_2 = \langle \sigma_c, [\sigma_a, \sigma_b]
\rangle$ and let $H = \langle \sigma_b, \sigma_c, [\sigma_a, \sigma_c],
[\sigma_a, \sigma_b] \rangle$.  Note that $D_1 \cong D_2 \cong D$ since
$\sigma_b$ and $[\sigma_a, \sigma_c]$ have order $2$ and $[\sigma_b, \,\,
[\sigma_a, \sigma_c]\,]$ is central of order $2$.  A similar argument holds
for $D_2$.  We have  $|H| = 32$, by observing the representation of
elements of $G_2$ at the end of the proof of Lemma 4.5.  We show
$\gal(L/F(\sqrt{a})) \cong D_1D_2 \cong DD$.  The groups $D_1$ and $D_2$
commute elementwise by Definition 4.4 and Lemma 4.5.  Therefore the
injections $D_1 \mapsto H$ and $D_2 \mapsto H$ induce a surjective
homomorphism $D_1 \times D_2 \mapsto H$.  The element \[ ([\sigma_b,
\,\,[\sigma_a, \sigma_c]\,], \,[\sigma_c, \,\,[\sigma_a, \sigma_b]\,] ) \]
lies in the kernel of the map $D_1 \times D_2 \mapsto H$.  Since
$[\sigma_b, \,\,[\sigma_a, \sigma_c]\,]$ and $[\sigma_c, \,\,[\sigma_a,
\sigma_b]\,]$ are the nontrivial elements in the center of $D_1$ and $D_2$,
it follows that $H$ is isomorphic to the central product $D_1D_2 \cong DD$.
It is clear that $H \subseteq \gal(L/F(\sqrt{a}))$ and therefore
$\gal(L/F(\sqrt{a})) = H \cong DD$.  This proves (2).

Suppose $E$ is a quadratic extension of $F$ such that $\gal(L/E) \cong DD$.
Then there is a subfield $M$ of $L$ such that $[L:M]=2$ and $\gal(M/E)
\cong ({\Bbb Z}/2{\Bbb Z})^4$.   (Note that $D \times D/ (\{1,g\} \times
\{1,g \}) \cong ({\Bbb Z}/2{\Bbb Z})^4$.)  Thus $M = E(\sqrt{\alpha_1}, \ldots,
\sqrt{\alpha_4})$, where $\alpha_i \in E$.  The Galois closure of
$E(\sqrt{\alpha_i})/F$ lies in $F^{(3)}$ for each $i$ and hence the Galois
closure of $M/F$ lies in $F^{(3)}$.  Thus $M$ lies in $F^{(3)}$ and it
follows $M=K$.  Thus $\gal(M/F) = \gal(K/F) \cong D \curlywedge D$.  Since
$D \curlywedge D$ contains a unique abelian subgroup of order $16$, it
follows $E = F(\sqrt{a})$ and this proves (3).

We have from before that $(a,b)_F = (a,c)_F =0$ and this implies $-a$ is
not rigid by Lemma 2.1.
\end{pf}

\section{Main Theorems}

\begin{Thm}
$F$ is not a rigid field if and only if either $G_1$ or $G_2$ occurs as a
Galois group over $F$.
\end{Thm}

\begin{pf}
First assume that $F$ is not a rigid field.  Then Lemma 2.1 implies there
exists $a \in \dot{F}\backslash F^2 $ such that $|D_F(\langle 1, -a
\rangle)/\dot{F}^2| \geq 4$.  Then either the statement in Theorem 3.7(4)
holds or the statement in Theorem 4.7(5) holds.  Therefore, either $G_1$ or
$G_2$ occurs as a Galois group over $F$.


Now assume that either $G_1$ or $G_2$ occurs as a Galois group over $F$.
Then by Theorems 3.7 and 4.7, there exists an element $a \in
\dot{F}\backslash F^2 $ such that $D_F(\langle 1, -a \rangle)$ contains at
least four square classes.   Then $F$ is nonrigid by Lemma 2.1.
\end{pf}

\begin{Thm}
Let $a \in \dot{F} \backslash (F^2 \cup -F^2)$.  Then the element $-a$ is
not rigid in $F$ if and only if at least one of the following two
statements holds.
\begin{enumerate}
\item An imbedding $F \subseteq F(\sqrt{a}) \subseteq E
\subseteq L$ exists such that $\gal(E/F) \cong C$ and $\gal(L/F) \cong G_1$.
\item An imbedding $F \subseteq F(\sqrt{a}) \subseteq L$ exists
such that $\gal(L/F(\sqrt{a})) \cong DD$ and $\gal(L/F) \cong G_2$.
\end{enumerate}
\end{Thm}

\begin{pf}
This follows from Proposition 3.8, Proposition 4.8, the construction in
Proposition 3.6 and the construction in Proposition 4.6.
\end{pf}

\begin{Cor}
Assume $F$ is nonrigid.  Let $v$ be a valuation on $F$ satisfying
$1+M \subseteq \dot{F}^2$ and $B(F) = \dot{F}^2 U$.
Then the following statements are equivalent.
\begin{enumerate}
\item $v$ is not $2$-divisible.
\item There exists an element $a \in \dot{F} \backslash \pm F^2$ such that
the two imbedding problems in Theorem 5.2 have no solutions for both
$F(\sqrt{a})$ and $F(\sqrt{-a})$.
\end{enumerate}
\end{Cor}

\begin{pf}
If $v$ is not $2$-divisible, then Corollary 2.3 implies there exists an
element $a \in \dot{F} \backslash \pm F^2$ such that $a$ is double rigid.
Now Theorem 5.2 implies that both statements in Theorem 5.2 fail
for both $F(\sqrt{a})$ and $F(\sqrt{-a})$.
Conversely, if both statements in Theorem 5.2 fail for both
$F(\sqrt{a})$ and $F(\sqrt{-a})$, then $a$ and $-a$ are rigid.
Thus $v$ is not $2$-divisible by Corollary 2.3.
\end{pf}

\section{The case $|\dot{F}/(\dot{F})^2| =4$}

\begin{Prop}
Assume $|\dot{F}/(\dot{F})^2| =4$.
\begin{enumerate}
\item The group $G_1$ occurs as a Galois group over $F$ if and only if $F$
is nonrigid.
\item If $F$ is nonreal, then $G_1$ occurs as a Galois group over $F$ if
and only if $Br_2(F) =0$ (i.e., $(a,b)_F=0$ for all $a,b \in F$).
\item If $F$ is formally real, then $G_1$ occurs as a Galois group over $F$
if and only if $F$ is uniquely ordered, which in this case is equivalent to
$F$ not being Pythagorean.
\end{enumerate}
\end{Prop}

\begin{pf}
(1).  If $G_1$ occurs as a Galois group over $F$, then $F$ is nonrigid by
Theorem 3.7.  Now assume $F$ is nonrigid.  Then Lemma 2.1 implies there
exists $a \in \dot{F}\backslash F^2 $ such that $|D_F(\langle 1, -a
\rangle)/\dot{F}^2| \geq 4$. Thus $D_F(\langle 1, -a \rangle) =\dot{F}$.
There exists  $b \in F$ such that $\dot{F} = \{1,a,b,ab\}\dot{F}^2$ and so
statement (3) of Theorem 3.7 holds and $G_1$ occurs as a
Galois group over $F$.

(2).  Let $F= F^2 \cup aF^2 \cup bF^2 \cup abF^2$.  Then $[F(\sqrt{a},
\sqrt{b}):F] =4$.  If $Br_2(F) =0$, then $(a,a)_F = (a,b)_F =0$.  Thus
$G_1$ occurs as a Galois group over $F$ by Theorem 3.7.

Now assume $G_1$ occurs as a Galois group over $F$.  Then there exist $a,b
\in F$ such that $[F(\sqrt{a}, \sqrt{b}):F] =4$ and $(a,a)_F = (a,b)_F =0$.
It follows that $F= F^2 \cup aF^2 \cup bF^2 \cup abF^2$ and thus $(a,
\dot{F})_F = 0$.  Suppose first that either $-1 \in F^2$ or $a \in -F^2$.
Then $(-1, \dot{F})_F =0$.  This implies $(b,b)_F = (-1,b)_F = 0$.  Then
$(b, \dot{F})_F =0$ (since $(b,a)_F=0$ and $(b,b)_F = 0$) and so $(ab,
\dot{F})_F =0$.  Therefore $Br_2(F) =0$.

Now assume $-1 \notin F^2$ and $a \notin -F^2$.  Then we can assume $b=-1$
(since we also have $(a,ab)_F=0$) and so $F=F^2 \cup -F^2 \cup aF^2 \cup
-aF^2$.  It is clear that $F^2 \cup aF^2 \subseteq D_F(\langle 1,a \rangle)$
and we have $F^2 \cup aF^2 \subseteq D_F(\langle 1,1 \rangle)$ since
$(-1,a)_F=0$.  Since $F^2 \cup aF^2 $ is not an ordering of $F$, we must
have either $F^2 \cup aF^2 \subsetneq D_F(\langle 1,a \rangle)$ or
$F^2 \cup aF^2 \subsetneq D_F(\langle 1,1 \rangle)$.  In the first case we
have $(-a, \dot{F})_F =0$ and in the second case we have $(-1, \dot{F})_F =0$.
Each case implies $Br_2(F) =0$.

(3). Since $F$ is formally real and $|\dot{F}/(\dot{F})^2| =4$, there is an
ordering such that $a$ is positive in this ordering and $F=F^2 \cup -F^2
\cup aF^2 \cup -aF^2$.  If $F$ is not Pythagorean then $a \in D_F(\langle
1,1 \rangle)$, and thus $F$ is uniquely ordered since $F^2 \cup aF^2$ is
contained in the positive cone of any ordering of $F$.  If $F$ is
Pythagorean, then one checks that $-1,a \notin D_F(\langle 1,-a \rangle)$
and so $F^2 \cup -aF^2$ is also the positive cone of an ordering.  Thus $F$
is not uniquely ordered.  The result now follows from Theorem 3.7.
\end{pf}

\begin{Prop}
Assume $|\dot{F}/(\dot{F})^2| =4$, $F$ is formally real and not
Pythagorean.  Then $[F^{(3)}:F] = 16$ with $\gal(F^{(3)}/F) \cong D
\curlywedge C$, and $[F^{\{3\}}:F] =32$ with $\gal(F^{\{3\}}/F) \cong G_1$.
\end{Prop}

\begin{pf}
We have $F = F^2 \cup -F^2 \cup aF^2 \cup -aF^2$ and $D_F(\langle 1,1
\rangle) = F^2 \cup aF^2$ where $F^2 \cup aF^2$ is the positive cone of the
unique ordering on $F$.  Let $E=F(\sqrt{-1})$.  Then the square class exact
sequence ([La], p. 202) applied to $E/F$ gives $|\dot{E}/(\dot{E})^2| =4$
since $D_F(\langle 1,1 \rangle)$ contains exactly two square classes.  Since
$F^{(2)}=E(\sqrt{a})$, the square class exact sequence applied to
$F^{(2)}/E$ gives $|\dot{F}^{(2)} /(\dot{F}^{(2)})^2| \leq
\frac{1}{2}|\dot{E}/(\dot{E})^2| =8$.  Therefore $[F^{\{3\}}:F^{(2)}] =
|\dot{F}^{(2)} /(\dot{F}^{(2)})^2| \leq 8$ and so $[F^{\{3\}}:F] \leq 32$.
Since $G_1$ occurs as a Galois group for some extension in $F^{\{3\}}$, it
follows $[F^{\{3\}}:F] =32$ and $\gal(F^{\{3\}}/F) \cong G_1$.  Since $D
\curlywedge C$ occurs as a Galois group over $F$ and $F^{(3)}
\subsetneq F^{\{3\}}$, it follows $\gal(F^{(3)}/F) \cong D \curlywedge C$.
\end{pf}

\begin{Prop}
Assume $|\dot{F}/(\dot{F})^2| =4$, $F$ is nonreal and $Br_2(F) =0$ (so that
$G_1$ occurs as a Galois group over $F$).  Then $[F^{(3)}:F] = 32$ and
$[F^{\{3\}}:F] = 128$.
\end{Prop}

\begin{pf}
We have $F^{(2)} = F(\sqrt{a}, \sqrt{b})$ for some $a,b \in F$.  Let
$E=F(\sqrt{a})$.  The square class exact sequence applied to $E/F$ gives
$|\dot{E}/(\dot{E})^2| = \frac{1}{2}|\dot{F}/(\dot{F})^2| =8$, since $u(F)
=2$ and so binary
quadratic forms over $F$ are universal.  Since $F$ is nonreal, it is known
that $u(E)=2$ also and thus the square class exact sequence applied to
$F^{(2)}/E$ gives $|\dot{F}^{(2)} /(\dot{F}^{(2)})^2| = \frac{1}{2} \cdot
8^2 = 32$.  Thus $[F^{\{3\}}:F^{(2)}] =
|\dot{F}^{(2)} /(\dot{F}^{(2)})^2| = 32$ and so $[F^{\{3\}}:F] =128$.

An argument in the proof of Proposition 3.6 shows that $\gal(F^{(3)}/F)$
is generated by any two automorphisms $\sigma_a, \sigma_b$ that satisfy
$\sigma_a(\sqrt{a}) = -\sqrt{a}$, $\sigma_a(\sqrt{b}) = \sqrt{b}$,
$\sigma_b(\sqrt{a}) = \sqrt{a}$, $\sigma_b(\sqrt{b}) = -\sqrt{b}$.
Then Corollary 2.8(2) implies $|\gal(F^{(3)}/F)| \leq 32$.
One can check that $C^aC^b \subseteq F^{(3)}$ and $D^{a,b} \subseteq
F^{(3)}$.  Since $[C^aC^b :F] =16$, $\gal(C^aC^b/F)$ is abelian and
$\gal(D^{a,b}/F)$ is nonabelian, it follows $[F^{(3)} :F] =32$.
\end{pf}

Generators and relations for $\gal(F^{\{3\}}/F)$ and
$\gal(F^{(3)}/F)$ can be found using Proposition 2.7 and Corollary 2.8.

We now briefly consider the case when $|\dot{F}/(\dot{F})^2| =4$ and $F$ is
a rigid field.  Then $F^{(3)} = F^{\{3\}}$ by [LS], although when
$|\dot{F}/(\dot{F})^2| =4$, this is easily obtained by direct calculation.
In the proposition below, $Q_8 \curlywedge C$ denotes the pullback of the
system of nontrivial homomorphisms $Q_8 \mapsto {\Bbb Z}/2{\Bbb Z}$ and $C
\mapsto {\Bbb Z}/2{\Bbb Z}$, where $Q_8$ is the quaternion group of order
$8$.  The group $Q_8 \curlywedge C$ is also the pullback of the system $D
\mapsto {\Bbb Z}/2{\Bbb Z}$ and $C \mapsto {\Bbb Z}/2{\Bbb Z}$, where the
kernel of $D \mapsto {\Bbb Z}/2{\Bbb Z}$ is the cyclic group of order $4$.
See [GSS] for more details.

\begin{Prop}
Assume $|\dot{F}/(\dot{F})^2| =4$ and $F$ is a rigid field.
\begin{enumerate}
\item If $F$ is formally real, then $[F^{\{3\}} :F]=8$ and
$\gal(F^{\{3\}}/F) \cong D$.
\item If $F$ is nonreal, then $[F^{\{3\}} :F]= 16$.
\begin{enumerate}
\item If $-1 \in F^2$, then $\gal(F^{\{3\}}/F) \cong C \times C$.
\item If $-1 \notin F^2$, then $\gal(F^{\{3\}}/F) \cong Q_8 \curlywedge C$.
\end{enumerate}
\end{enumerate}
\end{Prop}

\begin{pf}
(1). We have $F$ is Pythagorean by Proposition 6.1 and
$F^{(2)} = F(\sqrt{-1}, \sqrt{a})$ for some $a \in F$.  Two applications of
the square class exact sequence give
$|\dot{F(\sqrt{-1})}/\dot{F(\sqrt{-1})}^2| =2$ and
$|\dot{F^{(2)}}/(\dot{F^{(2)}})^2| =2$.  Thus $[F^{\{3\}} :F^{(2)}]= 2$ and
$[F^{\{3\}} :F] =8$.  We have $\gal(F^{\{3\}}/F) \cong D$ since $D^{a, -a}$
is a dihedral extension of $F$.

(2). When $F$ is nonreal, two applications of the square class exact
sequence give $|\dot{F^{(2)}}/(\dot{F^{(2)}})^2| = 4$ and thus $[F^{\{3\}}
:F] = 16$.  Part (a) follows since each quadratic extension of $F$ can be
imbedded in a cyclic quartic extension of $F$. One way to show (b) is to
see that $C^{-1}$ and $D^{a, -a}$ are extensions of $F$ and that
$\gal(D^{a, -a}/F(\sqrt{-1})) \cong C$.  Then $\gal(F^{\{3\}}/F) \cong Q_8
\curlywedge C$.
\end{pf}


\end{document}